\documentstyle[10pt,amscd]{amsart}

\newtheorem{proposition}{Proposition}

\newtheorem{remark}{Remark}

\begin{document}

\title{On the cobar construction of a bialgebra}

\author{Tornike Kadeishvili}


\keywords{Cobar construction, DG-bialgebra, Homotopy G-algebra, }

\begin{abstract}
We show that the cobar construction of a DG-bialgebra is a homotopy G-algebra. This implies that the bar construction of this cobar is a DG-bialgebra as well.
\end{abstract}

\maketitle


\section{Introduction}

The cobar construction $\Omega C$ of a  DG-coalgebra $(C,d:C\to C,\Delta:C\to C\otimes C)$ is, by  definition, a DG-algebra. Suppose now that $C$ is additionally equipped with a multiplication $\mu :C\otimes C\to C$ turning $(C,d,\Delta,\mu)$ into a  DG-bialgebra. How this multiplication reflects on  the cobar construction $\Omega C$?

It was shown by Adams \cite{Adams} that in mod 2 situation 
in this case the multiplication of  $\Omega C$ is homotopy commutative: there exists a $\smile_1$ product 
$$
\smile_1:\Omega C\otimes \Omega C\to \Omega C
$$
which satisfies the standard condition
\begin{equation}
\label{Steen}
\begin{array}{l}
d(a\smile_1b )=da\smile_1b + a\smile_1db + a\cdot  b + b \cdot  a,
\end{array}
\end{equation}
(since we work mod 2 the signes are ignored in whole paper).

In this note we show that this $\smile_1$ gives rise to a sequence of operations
$$
E_{1,k}:\Omega C\otimes (\Omega C)^{\otimes k}\to \Omega C,\ k=1,2,3,...
$$
which form on the cobar construction $\Omega C$ of a DG-bialgebra a structure of 
{\it homotopy  G-algebra} (hGa) in the sense of Gerstenhaber and Voronov \cite{GerVor}. 

There are two remarkable examples of homotopy G-algebras. 

The first one is the cochain complex of 1-reduced simplicial set $C^*(X)$. The operations $E_{1,k}$ here are dual to cooperations defined by Baues in \cite{Ba}, and the starting operation $E_{1,1}$ is  the classical Steenrod's $\smile_1$ product. 

The second example is the Hochschild cochain complex $C^*(U,U)$ of an associative algebra $U$. The operations $E_{1,k}$ here were defined in \cite{Kade} with the purpose to describe $A(\infty)$-algebras in terms of Hochschild cochains although the properties of those operations which where used as defining ones for the notion of homotopy G-algebra in \cite{GerVor} did not appear there. These operations where defined also in \cite{Getzler}. Again the starting operation $E_{1,1}$ is  the classical Gerstenhaber's circle product which is sort of $\smile_1$-product in the Hochschild complex. 

In this paper we present the third example of homotopy G-algebra: we construct the operations $E_{1,k}$ on the cobar construction 
$\Omega C$ of a DG-bialgebra $C$, and the starting operation $E_{1,1}$ is again classical, it is  Adams's $\smile_1$-product.
 
The notion of hGa was introduced in \cite{GerVor} as an additional structure on a DG-algebra $(A,d,\cdot)$ that induces a Gerstenhaber algebra structure on homology. The source of the  defining identities and the main example was Hochschild cochain complex $C^*(U,U)$. Another point of view is that hGa is a particular case of $B(\infty)$-algebra. This is an additional structure on a DG-algebra $(A,d,\cdot)$ that induces a DG-bialgebra structure on the bar construction $BA$.

We emphasize the third aspect of hGa: this is a structure which measures the noncommutativity of $A$. 
 There exists the classical tool which
measures the noncommutativity of a DG-algebra $(A,d,\cdot )$, namely the
Steenrod's $\smile_1$ product, satisfying the condition (\ref{Steen}).
The existence of such $\smile_1$ guarantees the commutativity of $H(A)$,
but $\smile _1$ product satisfying just the condition (\ref{Steen})
is too poor for most of applications. In many constructions some
deeper properties of $\smile_1$ are needed, for example the compatibility with the dot 
product of $A$ (the Hirsch formula)
\begin{equation}
\label{Hirsch}
(a\cdot b)\smile _1c + a\cdot (b\smile _1c)+(a\smile _1c)\cdot b=0.
\end{equation}
For a hGa $(A,d,\cdot,\{E_{1,k}\})$ the starting operation $E_{1,1}$ is a kind of $\smile_1$ product: it satisfies the conditions (\ref{Steen}) and (\ref{Hirsch}). As for the symmetric expression
$$
a\smile _1(b\cdot c)+ b\cdot (a\smile _1c)+(a\smile _1b)\cdot c,
$$
it is just {\it homotopical to zero} and the appropriate homotopy is the operation $E_{1,2}$. The defining conditions of a hGa which satisfy higher operations $E_{1,k}$ can  be regarded as  generalized Hirsch formulas. So we can say that a hGa is a DG-algebra with a "good" $\smile_1$ product.

\section{Notation and preliminaries}

We work over $Z_2$. 
For a graded $Z_2$-module $M$ we denote by $sM$ the suspension of $M$, i.e. 
$(sM)^i=M^{i-1}$ . Respectively $s^{-1}M$ denotes the desuspension of $M$, i.e. 
$(s^{-1}M)^i=M^{i+1}$.

A {\em differential
graded algebra } (DG-algebra) is a graded R-module $C=\{C^{i}\},\ i\in
 Z,$ with an associative multiplication $\mu :C^i\otimes
C^j\to  C^{i+j}$ and a homomorphism (a differential) $d:C^i
\to  C^{i+1}$ with $d^2=0$ and satisfying the Leibniz rule
$d(x\cdot y)=dx\cdot y+x\cdot dy,$ where $x\cdot y=\mu(x\otimes y)$. We
assume that a DG-algebra contains a unit $1\in C^0.$ A non-negatively
graded DG-algebra $C$ is {\em connected} if $C^0=Z_2.$ A connected  DG-algebra $C$
is {\em n-reduced}  if $C^i=0, 1\leq i\leq n.$ A DG-algebra is {\em
commutative} if $\mu=\mu T,$ where $T(x\otimes
y)=y\otimes x$.

 A {\em differential graded
coalgebra } (DG-coalgebra) is a graded $Z_2$-module $C=\{C_{i}\},\ i\in Z,$ with a coassociative comultiplication $\Delta :C\to  C\otimes
C$ and a homomorphism (a differential) $d:C_i \to  C_{i+1}$
with $d^2=0$ and satisfying  $\Delta d=(d\otimes id +id\otimes
d)\Delta.$ A DG-coalgebra $C$ is assumed to have a counit $\epsilon :C\to 
Z_2,\ \ (\epsilon \otimes id)\Delta =(id \otimes \epsilon)\Delta
=id. $ A non-negatively graded dgc $C$ is {\em connected} if
$C_0=Z_2.$ A connected  DG-coalgebra $C$ is {\em n-reduced} if $C_i=0, 1\leq
i\leq n.$ 

A {\em differential graded  bialgebra } (DG-bialgebra)
$(C,d,\mu , \Delta)$ is a DG-coalgebra $(C,d,\Delta)$ with a morphism of DG-coalgebras  
$\mu :C\otimes C\to  C$ turning $(C,d,\mu )$ into a DG-algebra.

\subsection{Cobar and Bar constructions}

Let $M$ be a graded $Z_2$-vector space with $M^{i\leq 0}=0$ and let $T(M)$ be the tensor algebra of $M$, i.e. 
$ T(M) =\oplus_{i=0}^{\infty}  M^{\otimes i}$. 
$T(M)$ is a free graded algebra: for a graded algebra $A$ and a homomorphism $\alpha:M\to A$ of degree zero there exists its {\it multiplicative extension}, a unique  morphism of graded algebras $f_{\alpha}:T(M)\to A$ such that $f_{\alpha}(a)=\alpha(a)$. The map $f_{\alpha}$ is given by 
$ f_{\alpha}(a_1\otimes...\otimes a_n)=\alpha(a_1)\cdot ... \cdot \alpha(a_n)$.

Dually, let $T^c(M)$ be the tensor coalgebra of $M$, i.e. 
$ T^c(M) =\oplus_{i=0}^{\infty}  M^{\otimes i}$, and the comultiplication 
$\nabla :  T^c(M)\to  T^c(M)\otimes  T^c(M)$ is given by
$$
\nabla (a_1\otimes...\otimes a_n)=\sum_{k=0}^n
(a_1\otimes...\otimes a_k)\otimes (a_{k+1}\otimes...\otimes a_n).
$$ 
$(T^c(M),\nabla)$ is a cofree graded coalgebra: for a graded coalgebra $C$ and a homomorphism $\beta:C\to M$ of degree zero there exists its {\it comultiplicative extension}, a unique  morphism of graded coalgebras $g_{\beta}:C\to T^c(M)$ such that $p_1g_{\beta}=\beta $, here $p_1:T^c(M)\to M$ is the clear projection. The map $g_{\beta}$ is given by 
$$
g_{\beta}(c)=\sum_n \beta(c^{(1)})\otimes...\otimes\beta(c^{(n)}),
$$
where $\Delta^n(c)=c^{(1)}\otimes ... \otimes c^{(n)}$ and $\Delta^n:C\to C^{\otimes n}$ is  $n$-th iteration of the diagonal $\Delta:C\to C\otimes C$, i.e. $\Delta^1=id,\ \Delta^2=\Delta,\ \Delta^n=(\Delta^{n-1} \otimes id)\Delta$.

Let $(C,d_C, \Delta )$ be a connected DG-coalgebra and $\Delta =id\otimes 1+1\otimes id+   {\Delta' }$. 
The (reduced) {\it cobar construction} $ \Omega C$ on $C$ is a DG-algebra whose underlying graded algebra is $T(sC^{>0})$. 
An element $(sc_1\otimes...\otimes sc_n)\in (sC)^{\otimes  n}\subset T(sC^{>0})$ is denoted by $[c_1,...,c_n]\in \Omega C$. The differential on $\Omega C$ is the sum $d= d_1+d_2$ which for 
a generator $[c]\in \Omega C $ is defined by
$d_1[c]=[d_C(c)] $
and
$d_2 [c]= \sum [c',c'']$ for $\Delta'(c)=\sum c'\otimes c''$, 
and extended   as a derivation.

Let $(A,d_A, \mu)$ be a 1-reduced  DG-algebra. The (reduced) {\it bar
construction} $ BA$ on $A$ is a DG-coalgebra whose underlying graded coalgebra is 
$T^c(s^{-1}A^{>0})$. Again an element 
$(s^{-1}a_1\otimes...\otimes s^{-1}a_n)\in (s^{-1}A)^{\otimes  n}\subset T^c(s^{-1}A^{>0})$ we denote as $[a_1,...,a_n]\in BA$. The differential of $BA$ is the sum $d=d_1 +d_2 $ which for  an element $[a_1,...a_n]\in BA$ is defined by 
$$
d_1[a_1,...,a_n]=\sum_{i=1}^{n} [a_1,...,d_Aa_i,...,a_n],
d_2 [a_1,...,a_n]= \sum_{i=1}^{n-1} [a_1,...,a_{i}\cdot a_{i+1},...,a_n].
$$


\subsection{Twisting cochains}\label{twist}

Let $(C,d,\Delta)$ be a dgc, $(A,d,\mu)$ a dga. 
A twisting cochain \cite{Brown} is a homomorphism $\tau:C\to A$ of degree +1 satisfying the Browns' condition 
\begin{equation}
\label{Brown}
d\tau+\tau d=\tau\smile \tau,
\end{equation}
where 
$\tau\smile \tau'=\mu_A(\tau\otimes\tau')\Delta $.
We denote by $T(C,A)$ the set of all twisting cochains $\tau:C\to A$.

There are universal twisting cochains $C\to \Omega C$ and $BA\to A$ being 
clear inclusion and projection respectively.
Here are essential consequences of the condition (\ref{Brown}):

\noindent (i)
The multiplicative extension 
$f_{\tau}:\Omega C\to A$ is a map of DG-algebras, so there is a bijection 
$T(C,A)\leftrightarrow Hom_{DG-Alg}(\Omega C,A)$;

\noindent (ii)
The comultiplicative extension 
$g_{\tau}:C\to BA$ is a map of DG0coalgebras, so there is a bijection 
$T(C,A)\leftrightarrow Hom_{DG-Coalg}(C,BA)$.


\section{Homotopy G-algebras}

\subsection{Products in the bar construction}

Let $(A,d,\cdot )$ be a 1-reduced DG-algebra and $BA$ it's bar construction.
We are interested in the structure of a multiplication
$$
\mu :BA\otimes BA\to BA,
$$
turning $BA$ into a DG-bialgebra, i.e. we require that

(i) $\mu $ is a DG-coalgebra map; 

(ii) is associative;

(iii) has the unit element $1_\Lambda \in \Lambda \subset BA$.

Because of the cofreeness of the tensor coalgebra
$BA=T^c(s^{-1}A)$,
a map of graded coalgebras
$$
\mu :BA\otimes BA\to BA
$$
is uniquely determined by the projection of degree +1
$$
E=pr\cdot \mu :BA\otimes BA\to BA\rightarrow A.
$$
Conversly, a homomorphism $E:BA\otimes BA\to A$ of degree +1 determines it's coextension, a graded coalgebra map $\mu _E:BA\otimes BA\to BA$ given by
$$
\mu _E=\sum_{k=0}^\infty (E\otimes ...\otimes E)\nabla _{BA\otimes BA}^k,
$$
where 
$\nabla _{BA\otimes BA}^k:BA\otimes BA\to (BA\otimes BA)^{\otimes k}$ 
is the k-fold iteration of the standard coproduct of tensor product of coalgebras 
$$
\nabla _{BA\otimes BA}=(id\otimes T\otimes id)(\nabla \otimes \nabla
):BA\otimes BA\rightarrow (BA\otimes BA)^{\otimes 2}.
$$

The map $\mu _E$ is a {\it chain map} (i.e. it is a map of DG-coalgebras) if and
only if $E$ is a twisting cochain in the sense of E. Brown, i.e. satisfies
the condition 
\begin{equation}
\label{brown}
dE+Ed_{BA\otimes BA}=E\smile E.
\end{equation}
Indeed, again because of the
cofreeness of the tensor coalgebra $BA=T^c(s^{-1}A)$ the condition 
$d_{BA}\mu _E=\mu _Ed_{BA\otimes BA}$ 
is satisfied if and only if it is
satisfied after the projection on $A$, i.e. if 
$pr\cdot d_{BA}\mu _E=pr\cdot \mu _Ed_{BA\otimes BA}$ 
but this condition is nothing else than the Brown's condition (\ref{brown}).

The same argument shows that the product $\mu _E$ {\it is associative} if and only if
$
pr\cdot \mu _E(\mu _E\otimes id)=pr\cdot \mu _E(id\otimes \mu _E),
$
or, having in mind $E=pr\cdot \mu _E$
\begin{equation}
\label{ass}
E(\mu _E\otimes id)=E(id\otimes \mu _E).
\end{equation}

A homomorphism $E:BA\otimes BA\rightarrow A$ consists of {\it components}
$$
\{\bar{E}_{p,q}:({s^{-1}A)^{\otimes p}}\otimes ({s^{-1}A)^{\otimes q}}\to A,\ p,q=0,1,2,...\} ,
$$
where $\bar{E}_{pq}$ is the restriction of $E$ on 
$({s^{-1}A)^{\otimes p}}\otimes ({s^{-1}A)^{\otimes q}}$. 
Each component $\bar{E}_{p,q}$ can be regarded as an operation 
$$
E_{p,q}:A^{\otimes p}\otimes A^{\otimes q}\to A,\ p,q=0,1,2,...\ .
$$
The value of $E_{p,q}$ on the element 
$(a_1\otimes ...\otimes a_p)\otimes (b_1\otimes ...\otimes b_q)$ 
we denote by 
$E_{p,q}(a_1,...,a_p;b_1...,b_q)$.

It is not hard to check that 
the multiplication $\mu _E$ induced by $E$ 
(or equivalently by a collection of multioperations $\{E_{p,q}\}$) {\it has the unit} $1_\Lambda \in \Lambda\subset BA$ if and only if
\begin{equation}
\label{unit}
E_{0,1}=E_{1,0}=id;\ \ E_{0,k}=E_{k,0}=0,\ k>1.
\end{equation}

So we can summarize:
\begin{proposition}
The multiplication $\mu _E$ induced by a collection of multioperations $\{E_{p,q}\}$ turns $BA$ into a DG-bialgebra, i.e. satisfies (i-iii), if and only if the conditions (\ref{brown}),
(\ref{ass}), and (\ref{unit}) are satisfied.
\end{proposition}


Let us interpret the condition (\ref{brown}) in terms of the components $E_{pq}$.

The restriction of (\ref{brown}) on $A\otimes A$ gives
\begin{equation}
\label{cup1}
\begin{array}{l}
dE_{1,1}(a;b)+E_{1,1}(da;b)+E_{1,1}(a;db)=
a\cdot b+b\cdot a.
\end{array}
\end{equation}
This condition coincides with the condition (\ref{Steen}), i.e. the operation $E_{1,1}$ is sort of $\smile _1$ product, which measures the noncommutativity of $A$. Below we denote $E_{1,1}(a;b)=a\smile _1b$.


The restriction on $A^{\otimes 2}\otimes A$ gives
\begin{equation}
\label{lhirsch}
\begin{array}{c}
d E_{2,1}(a,b;c)+E_{2,1}(da,b;c)+E_{2,1}(a,db;c)+E_{2,1}(a,b;dc)= \\
(a\cdot b)\smile _1c+a\cdot (b\smile _1c)+(a\smile _1c)\cdot b,
\end{array}
\end{equation}
this means, that this $\smile _1$ satisfies the {\it left
Hirsch formula} (\ref{Hirsch}) up to homotopy and the appropriate homotopy is the operation
$E_{2,1}$.


The restriction on $A\otimes A^{\otimes 2}$ gives:
\begin{equation}
\label{rhirsch}
\begin{array}{c}
dE_{1,2}(a;b,c)+E_{1,2}(da;b,c)+ E_{1,2}(a;db,c)+E_{1,2}(a;b,dc)= \\
a\smile _1(b\cdot c)+(a\smile _1b)\cdot c+ 
b\cdot (a\smile _1c),
\end{array}
\end{equation}
this means, that this $\smile _1$ satisfies the {\it right
Hirsch formula} (\ref{Hirsch}) up to homotopy and the appropriate homotopy is the operation
$E_{1,2}$.


Generally the restriction of (\ref{brown}) on $A^{\otimes m}\otimes A^{\otimes n}$ gives:
\begin{equation}
\label{Emn}
\begin{array}{c}
d E_{m,n}(a_1, ..., a_m;b_1, ..., b_n)+ 
\sum_iE_{m,n}(a_1, ..., d a_i, ..., a_m;b_1,.., b_n) \\
+\sum_iE_{m,n}(a_1, ..., a_m;b_1, ..., d b_i,.., b_n)= \\
a_1\cdot E_{m-1,n}(a_2, ..., a_m;b_1, ...,b_n)+
E_{m-1,n}(a_1, ..., a_{m-1};b_1, ..., b_n)\cdot a_m
\\
+b_1\cdot E_{m,n-1}(a_1, ..., a_m;b_2, ...,b_n)+
E_{m,n-1}(a_1, ..., a_m;b_1, ..., b_{n-1})\cdot b_m+\\
\sum_iE_{m-1,n}(a_1, ..., a_i\cdot a_{i+1}, ...,a_m;b_1, ..., b_n)+ \\
\sum_iE_{m,n-1}(a_1, ..., a_m;b_1, ..., b_i\cdot b_{i+1}, ..., b_n)+ \\
\sum_{p=1}^{m-1}\sum_{q=1}^{n-1} E_{p,q}(a_1, ..., a_p;b_1,.., b_q)\cdot 
E_{m-p,n-q}(a_{p+1}, ..., a_m;b_{q+1}, ..., b_n).
\end{array}
\end{equation}

Now let us interpret the associativity condition (\ref{ass}) in terms of  
the components $E_{p,q}$.
The restriction of (\ref{ass}) on $A\otimes A \otimes A$
gives
\begin{equation}
\label{cup1assoc}
\begin{array}{c}
(a\smile _1b)\smile  _1c + a\smile _1(b\smile _1c)=
E_{1,2}(a;b,c)+E_{1,2}(a;c,b)+ \\
E_{2,1}(a,b;c)+E_{2,1}(b,a;c).
\end{array}
\end{equation} 
Generally the restriction of (\ref{ass}) on 
$A^{\otimes k}\otimes A^{\otimes l} \otimes A^{\otimes m}$ gives
\begin{equation}
\label{Eassoc}
\begin{array}{c}
\sum_{r=1}^{l+m}\sum_{l_1+...+l_r=l, m_1+...+m_r=m}\\
E_{k,r}(a_1 ,  ... ,  a_k;E_{l_1,m_1}(b_1 ,  ... , 
b_{l_1};c_1 ,  ... ,  c_{m_1}) ,  ... ,  \\
E_{l_r,m_r}(b_{l_1+...+l_{r-1}+1} ,  ... , 
b_l;c_{m_1+...+m_{r-1}+1} ,  ... ,  c_m)= \\
\sum_{s+1}^{k+l}\sum_{k_1+...+k_s=k, l_1+...+l_s=l}\\
E_{s,m}(E_{k_1l_1}(a_1 ,  ... ,  a_{k_1};b_1 ,  ... , 
b_{l_1}) ,  ... ,  \\
E_{k_s,l_s}(a_{k_1+...+k_{s-1}+1} ,  ... , 
a_k;b_{l_1+...+l_{s-1}+1} ,  ... ,  b_l);c_1 ,  ... ,  c_m)
\end{array}
\end{equation}


We define a {\it Hirsch algebra} as a DG-algebra $(A,d,\cdot )$ endowed with a sequence of multioperations $\{E_{p,q}\}$ satisfying (\ref{unit}), (\ref{Emn}). 
This name is inspired by the fact that the defining condition (\ref{Emn}) 
can be regarded as generalizations of classical Hirsch formula (\ref{Hirsch}). This notion was used in \cite{KS}, \cite{KSan}.

A Hirsch algebra we call  {\it associative} if in addition the condition (\ref{Eassoc}) is satisfied.

This structure is a particular case of $B_\infty $-algebra, see below. Moreover turn the notion of {\it homotopy G-algebra}, described below, is a particular case of an associative Hirsch algebra.

\subsection{Some particular cases}

For a Hirsch algebra $(A,d,\cdot ,\{E_{p,q}\})$ the operation 
$E_{1,1}=\smile_1$ satisfies (\ref{Steen}), so this structure can be considered as a tool which measures  the noncommutativity of the product $a\cdot b$ of $A$. 
We distinguish various levels of ''noncommutativity'' of $A$ according to
the form of $\{E_{p,q}\}$.

\noindent{\bf Level 1.} Suppose for the collection $\{E_{p,q}\}$ all the operations except 
$E_{0,1}=id$ and $E_{1,0}=id$ are trivial. Then it follows from (\ref{cup1}) that 
in this case $A$ is a {\it strictly} commutative DG-algebra.

\noindent{\bf Level 2.} Suppose all operations except $E_{0,1}=id,\ E_{1,0}=id$ and $E_{1,1}$ are trivial. In this case $A$ is endowed with a ''strict'' $\smile _1$
product $a\smile _1b=E_{1,1}(a;b)$: 
the condition (\ref{Emn}) here degenerate to the following 4 conditions
$$
\begin{array}{l}
d(a\smile_1b )=da\smile_1b + a\smile_1db + 
a\cdot  b + b \cdot  a,
\end{array}
$$
$$
(a\cdot b)\smile _1c + a\cdot (b\smile _1c) + (a\smile _1c)\cdot b=0,
$$
$$
a\smile _1(b\cdot c)+ b\cdot (a\smile _1c) +(a\smile _1b)\cdot c=0,
$$
$$
(a\smile _1c)\cdot (b\smile _1d)=0.
$$
The condition (\ref{Eassoc}) degenerates to the associativity $\smile _1$
$$
a\smile _1(b\smile _1c)=(a\smile _1b)\smile _1c.
$$

As we see in this case we have very strong restrictions on $\smile _1$-product. An
example of DG-algebra with such strict $\smile _1$ product is 
$(H^{*}(SX,Z_2),d=0)$ with $a\smile _1b=0$ if $a\neq b$ and $a\smile
_1a=Sq^{|a|-1}a$, another example is $C^{*}(SX,CX)$,
where $SX$ is the suspension and $CX$ is the cone of a space $X$ (see \cite{San}).

\noindent{\bf Level 3.} Suppose all operations except 
$E_{0,1}=id,\ E_{1,0}=id$  and $E_{1,k},\ k=1,2,3,...$ are trivial.
In this case the condition (\ref{Emn}) degenerates into two conditions: at 
$A\otimes A^{\otimes k}$
\begin{equation}
\label{E1n}
\begin{array}{c}
dE_{1,k}(a;b_1, ..., b_k)+ E_{1,k}(da;b_1, ...,b_k)+
\sum_i 
E_{1,k}(a;b_1, ..., db_i, ..., b_k)= \\
b_1\cdot E_{1,k-1}(a;b_2, ..., b_k)
+ 
E_{1,k-1}(a;b_1, ..., b_{k-1})\cdot b_k+\\
\sum_iE_{1,k-1}(a;b_1, ...,
b_i\cdot b_{i+1}, ..., b_k),
\end{array}
\end{equation}
and at $A^{\otimes 2}\otimes A^{\otimes k}$
\begin{equation}
\label{E2n}
\begin{array}{c}
E_{1,k}(a_1\cdot a_2;b_1,.., b_k)= 
a_1\cdot E_{1,k}(a_2;b_1, ..., b_k)+ E_{1,k}(a_1;b_1, ..., b_k)\cdot a_2+\\
\sum_{p=1}^{k-1} 
E_{1,p}(a_1;b_1, ..., b_p)\cdot
E_{1,m-p}(a_2;b_{p+1}, ..., b_k);
\end{array}
\end{equation}
moreower at $A^{\otimes n>2}\otimes A^{\otimes k}$ the condition is trivial.
In particular the condition (\ref{lhirsch}) here degenerates to Hirsch formula 
(\ref{Hirsch}).

The associativity condition (\ref{Eassoc}) in this case looks as
\begin{equation}
\label{E1assoc}
\begin{array}{c}
E_{1,n}(E_{1,m}(a ;b _1, ..., b _m);c_1, ..., c _n)=
\sum_{0\leq i_1\leq ...\leq i_m\leq n}     \sum_{0\leq n_1+...+n_r\leq n}
\\
E_{1,n-(n_1+...+n_j)+j}(a ;c _1, ..., 
c _{i_1},E_{1,n_1}(b _1;c _{i_1+1}, ..., c _{i_1+n_1}), c _{i_1+n_1+1},
..., \\
c _{i_2},E_{1,n_2}(b _2;c _{i_2+1}, ..., c _{i_2+n_2}), c _{i_2+n_2+1},
...,\\
c _{i_m},E_{1,n_m}(b _m;c _{i_m+1}, ..., c _{i_m+n_m}), c _{i_m+n_m+1},...,c _n),
\end{array}
\end{equation}
In particular the condition (\ref{cup1assoc}) here degenerates to  
\begin{equation}
\label{cup1ass}
\begin{array}{c}
(a\smile _1b)\smile  _1c + a\smile _1(b\smile _1c)=
E_{1,2}(a;b,c)+E_{1,2}(a;c,b).
\end{array}
\end{equation}

The structure of this level coincides with the notion of 
{\it Homotopy G-algebra}, see below.

\noindent{\bf Level 4.} As the last level we consider a Hirsch algebra structure with no restrictions. An example of such structure is the cochain complex of a 1-reduced cubical set. Note that it is a {\it nonassociative} Hirsch algebra.

\subsection{$B_\infty $-algebra}

The notion of $B_\infty -$algebra was introduced in \cite{Ba}, \cite{GJ} as
an additional structure on a DG-algebra $(A,\cdot ,d)$ which turns the
tensor coalgebra $T^c(s^{-1}A)=BA$ into a DG-bialgebra. So it requires a
new differential
$$
\widetilde{d}:BA\rightarrow BA
$$
(which should be a coderivation with respect to standard coproduct of $BA$)
and a new associative multiplication
$$
\widetilde{\mu }:(BA,\widetilde{d})\otimes (BA,\widetilde{d})\rightarrow (BA,%
\widetilde{d})
$$
which should be a map of DG-coalgebras, with $1_\Lambda \in \Lambda \subset
BA$ as the unit element.

It is known that
such $\widetilde{d}$ specifies on $A$ a structure of $A_\infty $-algebra in
the sense of Stasheff \cite{Sta}, namely a sequence of operations 
$\{m_i:\otimes ^iA\rightarrow A,i=1,2,3,...\}$ subject of appropriate
conditions.

As for the new multiplication $\widetilde{\mu }$, it follows from the above
considerations, that it is induced by a sequence of operations $\{E_{pq}\}$
satisfying (\ref{unit}), (\ref{Eassoc}) and the modified condition (\ref{Emn}%
) with involved $A_\infty $-algebra structure $\{m_i\}$.

Thus the structure of associative Hirsch algebra is a particular $B_\infty $-algebra
structure on $A$ when the standard differential of the bar construction
$d_B:BA\rightarrow BA$ {\it does not change}, i.e. $\widetilde{d}=d_B$ (in
this case the corresponding $A_\infty $-algebra structure is degenerate:
$m_1=d_A,m_2=\cdot ,m_3=0,m_4=0,...$).

Let us mention, that a twisting cochain $E$
satisfying (\ref{unit}) and (\ref{brown}), (but not (\ref{ass}) i.e. the
induced product in the bar construction is not strictly associative),
was constructed
in \cite{Khel} for the singular cochain complex of a topological space
$C^{*}(X)$ using acyclic models. The condition
(\ref{unit})  determines this twisting cochain $E$
uniquely up to standard equivalence (homotopy) of twisting cochains in the sense of N. Berikashvili \cite{Ber}.


\subsection{Strong homotopy commutative algebras}

The notion of strong homotopy commutative algebra (shc-algebra), as a tool
for measuring of noncommutativity of DG-algebras, was used in many papers: \cite{Munk}, \cite{Thoma}, etc.

A shc-algebra is a DG-algebra $(A,d,\cdot )$ with a given
twisting cochain
$\Phi:B(A\otimes A) \rightarrow  A$
which satisfies apropriate up to homotopy conditions of associativity and
commutativity. Compare with the Hirsch algebra structure which is represented by a twisting cochain $E:BA\otimes BA\to A$. Standard contraction of $B(A\otimes A)$ to $BA\otimes BA$ allows to establish connection between these two notions. 


\subsection{DG-Lie algebra structure in a Hirsch algebra}

A structure of an associative Hirsch algebra on $A$ induces on the homology $H(A)$ a
structure of Gerstenhaber algebra (G-algebra) (see \cite{Gerst0}, \cite{GerVor}, \cite{Vor}) 
which is defined
as a commutative graded algebra $(H,\cdot )$ together with a Lie bracket of
degree -1
$$
[\ ,\ ]:H^p\otimes H^q\rightarrow H^{p+q-1}
$$
(i.e. a graded Lie algebra structure on the desuspension $s^{-1}H$) that is a
biderivation: $[a,b\cdot c]=[a,b]\cdot c+b\cdot [a,c]$.

The existence of this structure in the homology $H(A)$ is seen by the following argument.

Let $(A,d,\cdot ,\{E_{p,q}\})$ be an associative Hirsch algebra, then in the desuspension 
$s^{-1}A$ there appears a structure of DG-Lie algebra: although the $\smile
_1=E_{1,1}$ is not associative, the condition (\ref{cup1assoc}) implies the pre-Jacobi
identity
$$
a\smile _1(b\smile _1c)+(a\smile _1b)\smile _1c=
a\smile _1(c\smile _1b)+(a\smile _1c)\smile
_1b
$$
This condition guarantees that the commutator
$[a,b]=a\smile _1b +b\smile _1a$
satisfies the Jacobi identity. Besides the condition (\ref{cup1}) implies
that $[\ ,\ ]:A^p\otimes A^q\rightarrow A^{p+q-1}$ is a chain map.
Thus on $s^{-1}H(A)$ there appears the structure of graded Lie algebra. The
up to homotopy Hirsh formulae (\ref{lhirsch}) and (\ref{rhirsch}) imply that
the induced Lie bracket is a biderivation.

\subsection{Homotopy G-algebra}

An associative Hirsch algebra of level 3 in the literature is known as
{\it Homotopy G-algebra }.

A {\it Homotopy G-algebra} in \cite{GerVor} and \cite{Vor} is defined as a
DG-algebra $(A,d,\cdot )$ with a given sequence of multibraces 
$a\{a_1,...,a_k\}$ which, in our notation, we regard as a sequence of
operations
$$
E_{1,k}:A\otimes (\otimes ^kA)\rightarrow A,\quad k=0,1,2,3,...
$$
which, together with $E_{01}=id$ satisfies the conditions (\ref{unit}), (\ref
{E1n}), (\ref{E2n}) and (\ref{E1assoc}).

The name {\it Homotopy G-algebra} is motivated by the fact that this
structure induces on the homology $H(A)$ the structure of G-algebra 
(as we have seen in the previous section such a structure appears even on the
homology of an associative Hirsch algebra).

The conditions (\ref{E1n}), (\ref{E2n}), and (\ref{E1assoc}) in \cite{GerVor} are called {\it higher homotopies}, {\it distributivity} and {\it higher pre-Jacobi identities} respectively.
As we have seen the first two conditions mean that $E:BA\otimes BA\to A$ is a twisting cochain, or equivalently $\mu_{E}:BA\otimes BA\to BA$ is a chain map, and the third one means that this multiplication is associative.

\subsection{Operadic description}

Appropriate language to describe such huge sets of operations is the operadic language. Here we use {\it surjection operad} $\chi$ and {\it Barratt-Eccles operad} ${\cal E}$ which are most convenient $E_{\infty}$ operads. For definitions we refer to \cite{Fresse}.

The operations $E_{1,k}$ forming hGa have nice description in the {\it surjection operad}, see \cite{Mac}, \cite{Mac2}, \cite{Fresse}. Namely, to the dot product corresponds the element 
$(1,2)\in \chi_0(2)$, to $E_{1,1}=\smile_1$ product corresponds $(1,2,1)\in \chi_1(2)$, to the operation $E_{1,2}$ the element $(1,2,1,3)\in \chi_2(3)$, etc. Generally to the operation $E_{1,k}$ corresponds the element
\begin{equation}
\label{surj}
E_{1,k}=(1,2,1,3,...,1,k,1,k+1,1)\in \chi_k(k+1).
\end{equation}
We remark  here that the defining conditions of a hGa (\ref{E1n}), (\ref{E2n}), 
(\ref{E1assoc}) can be expressed in terms of operadic structure (differential, symmetric group action and composition product) and the elements (\ref{surj}) satisfy these conditions 
{\it already in the operad} $\chi$. This in particular implies that {\it any $\chi$-algebra is automatically a hGa}. Note that the elements (\ref{surj}) together with $(1,2)$ generate the suboperad $F_2\chi$ which is equivalent to the little square operad. This fact and a hGa structure on the Hochschild cochain complex $C^*(U,U)$ of  an algebra $U$ are used by many authors to prove so called Deligne conjecture about the action of the little square operad on $C^*(U,U)$.

Now look at the operations $E_{p,q}$ which define a structure of Hirsch algebra. 
They {\it can not live} in $\chi$: it is enough to mention that the Hirsch formula (\ref{Hirsch}), as a part of defining conditions of  hGa, is satisfied in $\chi$, but for a Hirsch algebra this condition is  satisfied up  to homotopy $E_{2,1}$, see (\ref{lhirsch}). We belive that $E_{p,q}$-s live in the Barratt-Eccles operad ${\cal E}$. In particular direct calculation shows that 
$$
\begin{array}{ll}
&E_{1,1}=((1,2),(2,1))\in {\cal E}_{1}( 2);\\
&E_{1,2}=(( {\bf 1},2,3),(2,{\bf 1},3),(2,3,{\bf 1}))\in {\cal E}_2( 3);\\ 
&E_{2,1}=((1,2, {\bf 3}),(1,{\bf 3},2),({\bf 3},1,2))\in {\cal E}_2( 3); \\
&E_{1,3}=(({\bf 1},2,3,4),(2,{\bf 1},3,4),(2,3,{\bf 1},4),(2,3,4,{\bf 1}))\in {\cal E}_3(4); \\
&E_{3,1}=((1,2,3,{\bf 4}),(1,2,{\bf 4},3),(1,{\bf 4},2,3),({\bf 4},1,2,3))\in {\cal E}_3(4);
\end{array}
$$
and in general 
$$
\begin{array}{l}
E_{1,k}=(({\bf 1},2,...,k+1),...,(2,3,...,i,{\bf 1},i+1,...,k+1),...,(2,3,...,k+1,{\bf 1}));\\
E_{k,1}=((1,2,...,{\bf k+1}),...,(1,2,...,i,{\bf k+1},i+1,...,k),...,({\bf k+1},1,2,...,k)).
\end{array}
$$
As for other $E_{p,q}$-s we can indicate just
$$
\begin{array}{ll}
E_{2,2}=&((1,2,3,4),(1,3,4,2),(3,1,4,2),(3,4,1,2))+\\
        &((1,2,3,4),(3,1,2,4),(3,1,4,2),(3,4,1,2))+\\
        &((1,2,3,4),(1,3,2,4),(1,3,4,2),(3,1,4,2))+\\
        &((1,2,3,4),(1,3,2,4),(3,1,2,4),(3,1,4,2)).
\end{array}
$$

We remark that the operadic {\it table reduction } map $TR:{\cal E}\to \chi$, see \cite{Fresse}, maps $E_{k>1,1}$ and $E_{2,2}$ to zero, and $E_{1,k}\in {\cal E}_{k}(k+1)$ to 
$E_{1,k}\in {\chi}_{k}(k+1)$.


\section{Adams $\smile_1$-product in the cobar construction of a bialgebra}

Here we present the Adams $\smile_1$-product 
$\smile_1: \Omega A\otimes \Omega A\to \Omega A$
on the cobar construction $\Omega A$ of a DG-bialgebra $(A,d,\Delta:A\to A\otimes A,\mu:A\otimes A\to A)$ (see \cite{Adams}). This will be the first step in the construction of  hGa structure on $\Omega A$.

This $\smile_1$ product satisfies the Steenrod condition (\ref{Steen}) and the Hirsch formula (\ref{Hirsch}).

First we define the $\smile_1$-product of two elements $x=[a],y=[b]\in \Omega A$ of length 1 
as $[a]\smile_1[b]=[a\cdot b]$. Extending this definition by (\ref{Hirsch}) we obtain 
$$
\begin{array}{c}
[a_1,a_2]\smile_1[b]=([a_1]\cdot [a_2])\smile_1[b]=
[ a_1]\cdot ([a_2]\smile_1[b])+([a_1]\smile_1 [b])\cdot[a_2]=\\
 \ [ a_1]\cdot [a_2\cdot b]+[a_1\cdot b]\cdot[a_2]=
[ a_1,a_2\cdot b]+[a_1\cdot b,a_2].
\end{array}
$$
Further iteration of this process gives 
$$
[a_1,...,a_n]\smile_1[b]=\sum_i [a_1,...,a_{i-1},a_i\cdot b,a_{i+1},...,a_n].
$$
Now let's define $[a]\smile_1[b_1,b_2]=[a^{(1)}\cdot b,a^{(2)}\cdot b]$ where 
$\Delta a=a^{(1)}\otimes a^{(2)}$ is the value of the diagonal $\Delta :A\to A\otimes A$ on $[a]$. Inspection shows that the condition (\ref{Steen}) for short elements
$$
d([a]\smile_1[b])=d[a]\smile_1[b]+[a]\smile_1d[b]+[a]\cdot[b]+[b]\cdot[a].
$$
is satisfied.

Generally we define the $\smile_1$ product of an element  
$x=[a]\in \Omega A$ of length 1 and an element $y=[b_1,...,b_n]\in \Omega A$ of arbitrary length by 
$$
[a]\smile_1 [b_1,...,b_n]=[a^{(1)}\cdot  b_1,...,a^{(n)}\cdot  b_n],
$$
here $\Delta ^n(a)=a^{(1)}\otimes ...\otimes a^{(n)}$ is the n-fold
iteration of the diagonal $\Delta :A\to A\otimes A$ and $a\cdot  b=\mu(a\otimes b)$ is the product in $A$.

Extending this definition for the elements of arbitrary lengths 
$[a_1,...,a_m]\smile_1[b_1,...,b_n]$
by the Hirsch formula (\ref{Hirsch}) we obtain the general formula 
\begin{equation}
\label{Adams}
[a_1,...,a_m]\smile_1 [b_1,...,b_n]=
\sum_k[a_1,...,a_{k-1},a_k^{(1)}\cdot  b_1,...,a_k^{(n)}\cdot  b_n,a_{k+1},...,a_m].
\end{equation}
Of course so defined $\smile_1$ satisfies the Hirsch  formula (\ref{Hirsch}) automatically. 
It remains to prove the
\begin{proposition}
\label{steenmn}
This $\smile_1$ satisfies Steenrod condition (\ref{Steen}) 
$$
\begin{array}{c}
d_{\Omega}([a_1,...,a_m]\smile_1 [b_1,...,b_n])=\\
d_{\Omega}[a_1,...,a_m]\smile_1 [b_1,...,b_n]+
 [a_1,...,a_m]\smile_1 d_{\Omega}[b_1,...,b_n]+\\
\ [a_1,...,a_m,b_1,...,b_n]+[b_1,...,b_n,a_1,...,a_m].
\end{array}
$$
\end{proposition}
\noindent {\bf Proof.}
Let us denote this condition by $Steen_{m,n}$.
The first  step consists in direct checking of the conditions $Steen_{1,m}$ by induction on $m$. Furthermore, 
assume that $Steen_{m,n}$ is satisfied. Let us check the condition $Steen_{m+1,n}$ for $[a,a_1,...,a_m]\smile_1[b_1,...,b_n]$. We denote $[a_1,...,a_m]=x,\ [b_1,...,b_m]=y$. Using the Hirsch formula (\ref{Hirsch}), $Steen_{m,n}$, and $Steen_{1,n}$ we obtain:
$$\begin{array}{c}
d([a,a_1,...,a_m]\smile_1[b_1,...,b_n])=\\
d(([a]\cdot  x )\smile_1y)=
d([a]\cdot  (x\smile_1y)+([a]\smile_1y)\cdot  x)=\\
=d[a]\cdot (x\smile_1y)+
[a]\cdot (dx\smile_1y+x\smile_1dy+x\cdot y+y\cdot x)+\\
(d[a]\smile_1y+[a]\smile_1dy+[a]\cdot y+y\cdot [a])\cdot x+
([a]\smile_1y)dx=\\
d[a]\cdot (x\smile_1y)+
[a]\cdot  (dx\smile_1y)+[a]\cdot  (x\smile_1dy)+
[a]\cdot x\cdot  y+[a]\cdot  y\cdot  x+\\
(d[a]\smile_1y\cdot )x+([a]\smile_1dy\cdot )x+[a]\cdot y\cdot x+
y\cdot [a]\cdot x+([a]\smile_1y)dx.
\end{array}
$$
Besides, using Hirsch (\ref{Hirsch}) formula we  obtain 
$$
\begin{array}{c}
d[a,a_1,...,a_m]\smile_1[b_1,...,b_n]=\\
d([a]\cdot x)\smile_1y=
(d[a]\cdot x)\smile_1y + ([a]\cdot  dx)\smile_1y=\\
d[a]\cdot (x\smile_1y)+(d[a]\smile_1y )\cdot x+ 
[a]\cdot  (dx\smile_1y) + ([a]\smile_1y)\cdot  dx
\end{array}
$$
and
$$
\begin{array}{c}
[a,a_1,...,a_m]\smile_1d[b_1,...,b_n]=\\
([a]\cdot x)\smile_1dy=
 [a]\cdot (x\smile_1dy) + ([a]\smile_1dy)\cdot x,
\end{array}
$$
now it is evident that $Steen_{m+1,n}$ is satisfied. This completes the proof.


\section{Homotopy G-algebra structure on the cobar construction of a bialgebra}

Below we present a sequence of operations 
$$
E_{1,k}:\Omega A\otimes (\Omega A)^{\otimes k}\to \Omega A,
$$
which extends the above described $E_{1,1}=\smile_1$ to a structure of homotopy  G-algebra
on the cobar construction of a DG-bialgebra . This means that $E_{1,k}$-s satisfy the conditions  (\ref{E1n}), (\ref{E2n}) and (\ref{E1assoc}).

For $x=[a]\in \Omega A$ of length 1, $y_i\in \Omega A$ and $k>1$ we define $E_{1,k}([a];y _1,...y _k)=0$ and extend for an arbitrary $x =[a_1,...a_n]$ by  (\ref{E2n}). This 
gives 
$$
E_{1,k}([a_1,...,a_n];y _1,...,y _k)=0$$
for $n<k$ and 
$$
E_{1,k}([a_1,...,a_k];y _1,...,y _k)=
[a_1\diamond y_1,...,a_k\diamond y_k],
$$
here we use the notation $a\diamond (b_1,...,b_s)=(a^{(1)}\cdot  b_1,...,a^{(s)}\cdot   b_s)$, so using this notation $[a]\smile_1[b_1,...,b_s]=[a\diamond (b_1,...,b_s)]$.
Further iteration by (\ref{E2n}) gives the general formula  
\begin{equation}
\label{Ecob}
\begin{array}{c}
E_{1,k}([a_1,...,a_n];y _1, ...,y _k)=\\
\sum 
[a_1,...,a_{i_1-1},a_{i_1}\diamond y _1,a_{i_1+1},..., 
a_{i_k-1},a_{i_k}\diamond y _k, a_{i_k+1},...,a_n], 
\end{array}
\end{equation}
where the summation is taken over all $1\leq i_1<...<i_k\leq n$. 

Of course so defined operations $E_{1,k}$ automatically satisfy the condition (\ref{E2n}). 
It remains to prove the
\begin{proposition}
\label{hga}
The operations $E_{1,k}$ satisfy the conditions (\ref{E1n}) and (\ref{E1assoc}). 
\end{proposition}
\noindent {\bf Proof.} The condition (\ref{E1n}) is trivial for $x=[a]$ of length 1 and $k>2$. For $x=[a]$ and $k=2$ this condition degenerates to
$$
\begin{array}{c}
E_{1,1}([a]  ;y _1\cdot  y _2)+
y _1\cdot  E_{1,1}([a]  ;y _2)-
E_{1,1}([a]  ;y _1)\cdot  y _2=0
\end{array}
$$
and this equality easily follows from the definition of $E_{1,1}=\smile_1$. For a long $x=[a_1,...,a_m]$ the condition (\ref{E1n}) can be checked by induction on the length $m$ of $x$ using the condition (\ref{E2n}).

Similarly, the condition (\ref{E1assoc}) is trivial for $x=[a]$ of length 1 unless the case $m=n=1$ and in this case this condition degenerates to
$$
\begin{array}{c}
E_{1,1}(E_{1,1}(x  ;y );z)=E_{1,1}(x;E_{1,1}(y );z))+
E_{1,2}(x;y,z)+E_{1,2}(x;z,y).
\end{array}
$$
This equality easily follows from the definition of $E_{1,1}=\smile_1$. For a long $x=[a_1,...,a_m]$ the condition (\ref{E1assoc}) can be checked by induction on the length $m$ of $x$ using the condition (\ref{E2n}).

\begin{remark}
For a DG-coalgebra $(A,d,\Delta:A\to A \otimes A)$ there is a standard DG-coalgebra map 
$g_A:A\to B\Omega A$ from $A$ to the bar of cobar of $A$.
This map is the coextension of the universal twisting cochain $\phi_A:A\to \Omega A$ defined by $\phi(a)=[a]$ and is a weak equivalence, i.e. it induces an isomorphism of homology.
Suppose $A$ is a DG-bialgebra. Then the constructed sequence of operations $E_{1,k}$ define a multiplication $\mu_E:B\Omega A\otimes B\Omega A\to B\Omega A$ on the bar construction $B\Omega A$ so that it becomes a DG-bialgebra. Direct inspection shows that 
$g_A:A\to B\Omega A$ is multiplicative, so it is a weak equivalence of DG-bialgebras. 
Dualizing this statement we obtain a weak equivalence of DG-bialgebras 
$\Omega BA \to A$ which can be considered as a free (as an algebra) resolution of a DG-bialgebra $A$.
\end{remark}

\vspace{20mm}

kade@@rmi.acnet.ge
     
A. Razmadze Mathematical Institute of the Georgian Academy of Sciences
 
M. Alexidze str. 1, Tbilisi, 380093, Georgia

\end{document}